\documentclass[leqno,12pt]{article}

\usepackage{times}

\usepackage{amsmath} 
\usepackage{amsthm}
\usepackage{amssymb}

\newtheorem{thm}{Theorem}[section]

\newtheorem{cor}[thm]{Corollary}
\theoremstyle{definition}
\newtheorem{df}[thm]{Definition}
\theoremstyle{definition}
\newtheorem{rem}[thm]{Remark}
\theoremstyle{plain}

\makeatletter

\def\address#1#2{\begingroup
\noindent\parbox[t]{7.8cm}{%
\small{\scshape\ignorespaces#1}\par\vskip1ex
\noindent\small{\itshape E-mail address}%
\/: #2\par\vskip4ex}\hfill%
\endgroup}%
\makeatother

\title{\bf Steinness of the Fatou set for a rational map of the complex projective plane}
\author{\sc Kazutoshi Maegawa}
\date{}

\begin{document}

\maketitle


\begin{abstract}
For a dominant algebraically stable rational self-map of the complex projective plane of degree at least 2,\ 
we will consider three different definitions of Fatou set and show the equivalence of them. Consequently,\ it follows that all Fatou components are Stein. 
This is an improvement of an early result by Forn\ae ss and Sibony ([FS]).
\end{abstract}

\section{Introduction}

The study of complex dynamics in higher dimensions has been developped intensively,\ but some early papers such as [FS] are still helpful for recent studies.
In those papers,\ the dynamics of rational self-maps of the complex projective spaces are studied and 
several fundamental theorems concerning Fatou set and Julia set are contained. 

In this paper,\ we will improve two theorems which are contained in [FS] (Therem 5.2 and Theorem 5.7 in [FS]).
 When we compare rational maps of $\Bbb P^2$ with those of $\Bbb P^1$,\ one of remarkable 
differences is that those of $\Bbb P^2$ possibly have indeterminacy points. 
If we want to generalize the Fatou-Julia theory to two dimension,\ 
we have to find a suitable way to define Fatou set and Julia set paying attention to 
indeterminacy points. 
So,\ in this paper,\ for an algebraically stable rational self-map $f$ of $\Bbb P^2$ of degree at least 2,\ we will consider several different definitions of Fatou set and show the equivalence of them. 
One of them is the standard,\ that is,\ it is based on local equicontinuity of the iterates of $f$,\ 
i.e. Lyapunov stability. The others are more complex analytical. 
They are based on notions of normal family for meromorphic maps.
Throughout this paper,\ we often use Ivashkovich's theorems in [I] and 
Forn\ae ss-Sibony's theorems concerning Green currents in [S]. Particularly,\ by using an Ivashkovich's result about domains of normality for 
families of meromorphic maps,\ we can show that all Fatou components for $f$ are Stein manifolds. 
    
\section{Preliminaries}

We denote by $\Bbb P^k$ the complex projective space of complex dimension $k \ge 1$
 and we equip the Fubini-Study distance with $\Bbb P^k$.
By definition,\ there is a holomorphic map 
$$\pi : \Bbb C^{k+1} \setminus \{O\} \rightarrow \Bbb P^k $$
\noindent
which is a $\Bbb C^{*}$-bundle over $\Bbb P^k$ such that 
the fiber $\pi^{-1}(p)$ for $p \in \Bbb P^k$ is $L \setminus \{O\}$ 
where $L$ is a complex line in $\Bbb C^{k+1}$ through the origin $O$.
By definition,\ a rational self-map $f$ of $\Bbb P^k$ is lifted by $\pi$ to a polynomial self-map $F$ of $\Bbb C^{k+1}$
 which is of the form
$$F=(P_0,\cdots,P_{k})$$
\noindent
where $P_i\ \ 0 \le i \le k$ are homogeneus polynomials which have the same degree and have no common factors.
The degree ${\rm deg}(f)$ of $f$ is defined to be the degree of $P_i$.
A point $p \in \Bbb P^k$ such that $F(\pi^{-1}(p))=O$ is an indeterminacy point for $f$.
We denote by $I=I(f)$ the set of indeterminacy points for $f$. 
The dimension of $I$ is at most $k-2$ if $k \ge 2$. 
(If $k=1$,\ then $I$ is empty.)
We obtain the lift of the $n$-th iterate $f^{n}=f \circ \cdots \circ f$ ($n$ times)
by the cancellation of common factors of the component functions
 for $F^n$. 

\begin{df}([S])\ \ We say that $f$ is {\it algebraically stable (AS)} 
if $f^n$ maps no complex hypersurface in $\Bbb P^k$ to $I(f)$,\ for all $n\ge1$.
\noindent
This is equivalent to that ${\rm deg}(f^n)=({\deg}(f))^{n}$ for all $n\ge1$.  
\end{df} 

Let $\omega$ be the normalized Fubini-Study $(1,1)$ form in $\Bbb P^k$ and 
$f$ be an AS rational self-map $f$ of $\Bbb P^k$ of degree $d \ge 2$ which is dominant,\ i.e. $f(\Bbb P^k)=\Bbb P^k$.
Then,\ 
$$\frac{1}{d^{n}} (f^n)^{*}\omega \rightarrow T$$
\noindent
as $n \rightarrow \infty$ in the sense of currents,\ where $T$ is a positive closed (1,1) current
 such that $f^{*}(T)=dT$ ([S]).
We call $T$ the {\it Green (1,1) current} for $f$. 
The current $T$ is of the form
$$T=\omega+{\rm dd^{c}}v$$
\noindent
where $v$ is an integrable function in $\Bbb P^k$. 

It is obvious that for any $n \ge 1$,\ the map $f^n$ is holomorphic in the complement of the closure of $\bigcup_{n\ge 1} I(f^n)$.

\begin{df}
We define the {\it Fatou set} ${\cal F}$ to be the maximal open subset of $\Bbb P^k \setminus \overline{\bigcup_{n\ge 1} I(f^n)}$
 in which $\{f^n\}$ is locally equicontinuous. A connected component of ${\cal F}$ is called a {\it Fatou component}.
The complement ${\cal J}$ of ${\cal F}$ is called the {\it Julia set} for $f$. 
\end{df} 

The support of $T$ is extremely related to ${\cal J}$. In particular,\ the following equality is known.  
 
\begin{thm} ([FS],\ [U]) \ \ 
In case when $f$ is holomorphic,\ ${\cal J}={\rm supp}(T)$.
\end{thm} 
 
\section{Results}

We will show that for any dominant AS rational self-map $f$ of $\Bbb P^2$ of degree at least $2$,\ 
each Fatou component for $f$ is Stein. Let us begin with recalling some basic notions in 
function theory of several complex variables. 

\

We set two subsets $\Delta$ and $H$ in $\Bbb C^k$
as 
$$\Delta:=\{|z_1|,\cdots,|z_k|<1 \},$$
$$H:=\{|z_1|<r_1,\cdots,|z_{k-1}| <r_1,\ |z_k|<1 \}$$
$$\cup \{|z_1|<1,\cdots,\ |z_{k-1}| <1,\ r_2<|z_k|<1\}$$
\noindent
where $0<r_1,r_2<1$. 

\begin{df}
Let $U$ be a domain in $\Bbb P^k$.
We say that $U$ is {\it pseudoconvex} if
$h:\Delta \rightarrow \Bbb P^k$ is an injective holomorphic map
and $h(H) \subset U$, then $h(\Delta) \subset U$.

\end{df}  

\begin{df}
Let $U$ be a domain in $\Bbb P^k$. 
We say that $U$ is {\it Stein} if $U$ is holomorphically convex and holomorphically separable.
\end{df}
\noindent
Above,\ 'holomorphically convex' means that for any compact set in $U$,\ the convex hull by holomorphic functions in $U$ 
is again compact. 'Holomorphically separable' means that for any points $p,\ q \in U$,\ there exists a holomorphic 
function $h$ in $U$ such that $h(p) \neq h(q)$.

By Hartog's extention theorem of holomorphic functions,\ 
it can be shown that Stein domains are pseudoconvex. 
The converse is known as the Levi problem (or the Hartogs converse problem).
The following theorem due to Takeuchi is the solution for this problem. 
(Concerning the original work in case of domains in $\Bbb C^k$ by Oka, see [O].)  

\begin{thm} \label{Leviprob} ([T])\ \ 
Any pseudoconvex domain in $\Bbb P^k$,\ which is not $\Bbb P^k$,\ is Stein. 
\end{thm}

It is known that any Stein domain is biholomorphic to a closed complex submanifold in $\Bbb C^N$ for some $N$.
Although pseudoconvexity is just a condition on the shape of a domain,\ the theorem above gives 
plenty of complex analytical information of the domain.   

\

Concerning whether Fatou components are Stein or not,\ some partial answers are known.
In case of holomorphic maps,\ the affirmative answer was obtained in [FS] and [U]. 
In case of rational maps with indeterminacy points,\ Forn\ae ss-Sibony obtained 
the affirmative answer for some special maps. We need the following definition to explain this.

\begin{df}
Let $f$ be a dominant AS rational self-map of $\Bbb P^k$ of degree at least 2.
We say that a point $p \in \Bbb P^k$ is {\it regular} for $f$ if there exist a neighborhood $V$ of $p$ and a neighborhood $W$ of $I(f)$
 such that the orbit $\{f^{n}(V)\}_{n \ge 0}$ is disjoint from $W$.
\end{df}

\begin{rem}
A regular point is the same as a normal point in the sense of Sibony [S].
In this paper,\ we use 'regular' to avoid a confusion with normal points in the sense of Montel.
\end{rem} 

Forn\ae ss-Sibony showed that if all points in $\Bbb P^k \setminus \overline{\bigcup_{n\ge 1} I(f^n)}$ are regular for $f$,\ 
then all Fatou components for $f$ are Stein ([FS]). 

More recently,\ this result was improved. In [M],\ I established a dichotomy of Fatou components 
from a viewpoint of a dynamical relation with the indeterminacy set.

\begin{thm} ([M]) \ \ 
Let $f$ be a dominant AS rational self-map of $\Bbb P^k$ of degree at least 2.
Let $U$ be any Fatou component for $f$.
Then,\ either all points in $U$ are regular or no points in $U$ are regular. 
\end{thm}

We say that a Fatou component $U$ is {\it regular} if all points in $U$ are regular. 
To some extent,\ regular Fatou components are similar to Fatou components for holomorphic maps.
So,\ the following result can be obtained without much effort.

\begin{thm} ([M]) \ \ 
All regular Fatou components are Stein.
\end{thm} 

\

Here we will show our main theorem which is a complete answer to the question. 

\begin{thm} \label{main2}
Let $f$ be a dominant AS rational self-map of $\Bbb P^2$ of degree at least 2.
Then,\ any Fatou component for $f$ is Stein.
\end{thm}

Let us prepare some tools. 

\begin{df} ([I])
Let $X,Y$ be complex manifolds.
Let $\{g_n\}_{n \ge 1}$ be a sequence of meromorphic maps from $X$ to $Y$.
Let $\Gamma_n \subset X\times Y$ denote the graph of $g_n$. 
Let $g: X \rightarrow Y$ be a meromorphic map and $\Gamma\subset X \times Y$ be the graph of $g$. 

\begin{itemize}
\item[(i)] We say that $\{g_n\}_{n \ge 1}$ {\it strongly converges} to $g$ in $X$ if for any compact set $K \subset X$
$$\lim_{n \rightarrow \infty} \Gamma_n \cap (K \times Y)=\Gamma \cap (K \times Y)$$
\noindent
with respect to the Hausdorff metric.  
\item[(ii)] We say that $\{g_n\}_{n \ge 1}$ {\it weakly converges} to $g$ in $X$ if there is an analytic subset $A \subset X$
 of ${\rm codim}_{\Bbb C} A \ge 2$ such that $\{g_n\}_{n \ge 1}$ strongly converges to $g$ in $X \setminus A$.  
\end{itemize} 

\begin{rem}
The definition of strong convergence above is slightly modified from the original one in [I]. It is just a change in naming.
\end{rem} 

By using these two notions of convergence,\  
we can introduce notions of normality for a sequence of meromorphic maps in strong and weak senses,\ that is,\ we say that a sequence $\{f_n\}$ of meromorphic maps from $X$ to $Y$ is 
{\it strongly (resp. weakly) normal} if for any subsequence $\{f_{n_j}\}$ of $\{f_n\}$,\ there is a subsequence
 of $\{f_{n_j}\}$ which converges in $X$ in strong (resp. weak) sense.
Thus,\ for the iterates of a rational self-map of $\Bbb P^2$,\ 
we can define the {\it strong (resp. weak) Fatou set} ${\cal F}_s$ 
(resp. ${\cal F}_w$)
 as the maximal open subset of $\Bbb P^2$ in which the iterates is strongly (resp. weakly) normal.   
\end{df}

The following two theorems are contained in an Ivashkovich's paper [I].

\begin{thm} ([I]) \label{Iv}
Let $f$ be a rational self-map of $\Bbb P^2$. Then,\ the following (i) and (ii) hold.
\begin{itemize}
\item[(i)] ${\cal F}_{w}$ is pseudoconvex;
\item[(ii)] If ${\cal F}_{s} \neq {\cal F}_{w}$,\ then  ${\cal F}_{w} \supset \Bbb P^2 \setminus C$,\ where $C$ is a rational curve in $\Bbb P^2$.
\end{itemize}
\end{thm}

\begin{rem}
In [I],\ the assertion (i) is shown in a more general situation. That is,\ it holds for any sequence of meromorphic maps from any complex manifold to any compact
 K\"ahler manifold.
\end{rem} 

\begin{thm} ([I]) (Rouch\'e Principle) \label{RP} \ \ 
Let $X,Y$ be complex manifolds.  
Let $f_n :X \rightarrow Y,\ n\ge1,\ $ be meromorphic maps which strongly converge 
 in $X$ to a meromorphic map $f:X \rightarrow Y$ as $n \rightarrow \infty$.
Then,\ 
\begin{itemize}
\item[(a)]\ If $f$ is holomorphic,\ then for any relatively compact open subset $D$ in $X$,\ 
the restrictions $f_{n}|D$ for sufficiently large $n$ are holomorphic in $D$,\ and $\{f_n\}$ converges to $f$ uniformly on compact sets in $X$.
\item[(b)]\ If $\{f_n\}$ are holomorphic,\ then $f$ is also holomorphic,\ and $\{f_n\}$ converges uniformly on compact sets
 in $X$. 
\end{itemize} 
\end{thm}

We also need the following two theorems which come from the study of dynamics using pluripotential theory.
Especially,\ Theorem \ref{charge} plays a crucial role.

\begin{thm} {([S])} \label{Green}
Let $T$ be the Green (1,1) current for a dominant AS rational self-map $f$ of $\Bbb P^k$ of degree at least 2. Then,\ ${\cal F} \subset \Bbb P^k \setminus {\rm supp}(T)$.
\end{thm}

\begin{rem} \label{M}
We can easily obtain a slightly modified version of this theorem. That is,\ if a subsequence $\{f^{n_j}\}_{j\ge1}$ converges uniformly in 
an open set $U$ in $\Bbb P^k$,\ then $U \subset \Bbb P^k \setminus {\rm supp}(T)$.
\end{rem} 

Let $K$ be a closed set in $\Bbb P^k$ 
and $S$ be a positive current in $\Bbb P^k \setminus K$ with locally bounded mass near $K$.
The trivial extension of $S$ to $\Bbb P^k$ is obtained by 
setting the coefficients of $S$ to be $0$ on $K$.
(Note that the coefficients of $S$ are complex measures.)
  
\begin{thm} {([S])} \label{charge}
Let $T$ be the Green (1,1) current for a dominant AS rational self-map $f$ of $\Bbb P^k$ of degree at least 2. Then,\ $T$ does not charge on any complex hypersurface $V$,\ i.e. $T$ is equal to the trivial extension of $T|_{\Bbb P^k \setminus V}$.   
\end{thm}

\

\noindent
{\sc Proof of Theorem \ref{main2}:}\ By definition,\ it is easily shown that 
$${\cal F} \subset {\cal F}_{s} \subset {\cal F}_{w}.$$
Let us show that ${\cal F}_s = {\cal F}_w$.
Let $U_w$ be a connected component of ${\cal F}_{w}$.
Suppose that a sequence $\{f^{n_k}\}_{k\ge1}$ weakly converges in $U_w$.
Then,\ there is a discrete set $A \subset U_w$ such that $\{f^{n_k}\}_{k \ge 1}$ strongly converges to a meromorphic 
map $\psi$ in $U_w \setminus A$. Let ${I(\psi)}$ denote the indeterminacy set of $\psi$.
By definition,\ $\psi$ is holomorphic in $U_w \setminus {I(\psi)}$. 
Suppose that there is $m \ge 0$ such that $f^m$ is not holomorphic at $p \in U_w \setminus ({I(\psi)} \cup A)$.
Since $f$ is AS,\ it follows that $f^n$ is not holomorphic at $p$ for $\forall n \ge m$.
By Rouch\'e principle (Theorem \ref{RP}),\ $\psi$ is also not holomorphic at $p$.
This is a contradiction.
So,\ $f^n$ are holomorphic in $U_w \setminus ({I(\psi)} \cup A)$ for all $n \ge 1$.
Further,\ $\{f^{n_k}\}_{k \ge 1}$ converges locally uniformly to $\psi$ in $U_w \setminus ({I(\psi)} \cup A)$.
By Remark \ref{M},\ 
it follows that $U_w \setminus ({I(\psi)} \cup A)$ is contained in $\Bbb P^2 \setminus {\rm supp}(T)$.
Since ${I(\psi)} \cup A$ is discrete,\ the local potential function for $T$ should be pluriharmonic in $U_w$. 
Hence,\ it follows that $U_w \subset \Bbb P^2 \setminus {\rm supp}(T)$. Thus,\ ${\cal F}_{w} \subset \Bbb P^2 \setminus {\rm supp}(T)$. 
Suppose that ${\cal F}_s \neq {\cal F}_w$.
By (ii) in Theorem \ref{Iv},\ ${\cal F}_w \supset \Bbb P^2 \setminus C$,\ where $C$ is a rational curve.
Hence,\ $C \supset {\rm supp}(T)$.
However this is impossible because $T$ does not charge on $C$ (Theorem \ref{charge}).
Thus,\ it follows that ${\cal F}_s = {\cal F}_w$.  
 
Let us show that ${\cal F} = {\cal F}_s$.
Let $U_s$ be a connected component of ${\cal F}_s$. 
Let $\{f^{n_k}\}_{k \ge 1}$ be any subsequence of $\{f^n\}_{n \ge 1}$ which strongly converges in $U_s$. 
Denote the limit map by $\phi$ and let ${I(\phi)}$ be the indeterminacy set for $\phi$.
From Rouch\'e principle and $f$ being AS,\ in the same way as we have done above,\ it follows that 
$f^{n_k}$\ $k \ge 1$ are holomorphic in $U_s \setminus {I(\phi)}$ and 
the convergence is locally uniform. Hence,\ $U_s \setminus {I(\phi)} \subset \Bbb P^2 \setminus 
{\rm supp}(T)$. Since ${I(\phi)}$ is discrete,\ it follows that $U_s \subset \Bbb P^2 \setminus {\rm supp}(T)$.  
Suppose that ${I(\phi)}$ is nonempty and $q \in {I(\phi)}$. By Rouch\'e principle,\ there is an integer $m \ge 0$
such that $f^{m}(q) \in {I(f)}$. So,\ it follows that $q \in {\rm supp}(T)$ since the Lelong number for $T$ at $q$
 must be strictly positive (this can be shown by the assumption that $f$ is AS). This is a contradiction to that $U_s \subset \Bbb P^2 \setminus {\rm supp}(T)$.      
Hence,\ ${I(\phi)} = \emptyset$. Again by Rouch\'e principle,\ 
it follows that $\{f^{n_k}\}_{k \ge 1}$ are holomorphic in $U_s$ and converges locally uniformly to $\phi$ in $U_s$.
This leads to that $U_s \subset {\cal F}$. Hence,\ ${\cal F}_s \subset {\cal F}$.

Thus,
$${\cal F} = {\cal F}_{s} = {\cal F}_{w}$$
\noindent
are verified. By (i) in Theorem \ref{Iv} and these equalities,\ ${\cal F}$ is pseudoconvex. 
So,\ by Theorem \ref{Leviprob},\ it follows that ${\cal F}$ is Stein.

\qed

Let us state an important fact which we have already shown in the proof above.
 
\begin{thm} \label{Feq}
Let $f$ be a dominant AS rational self-map of $\Bbb P^2$ of degree at least 2. Then,\ 
$${\cal F} = {\cal F}_{s} = {\cal F}_{w}.$$
\end{thm}

\begin{rem}
It is known that ${\cal F} \subset \Bbb P^2 \setminus {\rm supp}(T)$ (Theorem \ref{Green}).
However,\ there is an example for which ${\cal F} \neq \Bbb P^2 \setminus {\rm supp}(T)$. See Example 2 in [RS]. 
\end{rem}

As a consequence of Theorem \ref{main2},\ the connectivity of the Julia set follows. 
It is a common property of pseudoconcave sets in the complex projective space.
    
\begin{cor}
Let $f$ be a dominant AS rational self-map of $\Bbb P^2$ of degree at least 2. Then,\ the Julia set ${\cal J}$ is connected.
\end{cor}
\noindent
{\sc Proof:}\ \ The same way as in Theorem 5.2 in [FS] is valid. \qed

\

Let us note that Theorem \ref{Feq} does not necessarily hold 
in case of non AS maps. An example is a birational map $f([x:y:z])=[yz:zx:xy]$. 
When $m$ is even,\ $f^m$ is the identy map,\ and when $m$ is odd,\ $f^m=f$. 
Hence,\ it follows immediately that ${\cal F}_{s} = {\cal F}_{w}=\Bbb P^2$.
However,\ since $f$ has indeterminacy points $[1:0:0],[0:1:0],[0:0:1]$,\ the set ${\cal F}$ 
is smaller than $\Bbb P^2$. In fact,\ ${\cal F}=\Bbb P^2 \setminus \pi(\{xyz=0\})$.

\bigskip
\bigskip
\bigskip

\address{
Institute for the Promotion of Excellence in Higher Education \\
Kyoto University  \\
Yoshida Nihonmatsu-cho, Sakyo-ku \\ 
Kyoto,\ 606-8501 \\
Japan \\ 
}
{km@math.h.kyoto-u.ac.jp}

\end{document}